\documentclass[12pt]{article}
\usepackage[strict]{changepage}
\usepackage{secdot}
\usepackage{sectsty}
\usepackage{extarrows,color}
\usepackage{authblk}
\usepackage{verbatim}
\allowdisplaybreaks[3]
\usepackage{bbm}
\usepackage[top=1.25in, bottom=1.25in, left=1.25in, right=1.25in]{geometry}
\usepackage{mdwlist}
\usepackage{amsmath,amsmath,graphicx}

\usepackage[symbol]{footmisc}
\DefineFNsymbols{nu}{* 1 2 3}\setfnsymbol{nu}   
\title{Statistical Estimation and Nonlinear Filtering in Environmental Pollution  \thanks{This research is supported by Southern University of Science and Technology Start up fund Y01286120 and NSFC grants 61873325 and 11831010. Corresponding author: Qizhu Liang}}
\author{Qizhu Liang \thanks{ (yb47422@connect.um.edu.mo).} , Jie Xiong  
 \thanks{Department of Mathematics, Southern University of Science and Technology, Shenzhen, China (xiongj@sustech.edu.cn).}  ~ and ~ Xingqiu Zhao \thanks{Department of Applied Mathematics, Hong Kong, China (xingqiu.zhao@polyu.edu.hk).} 
}

\date{}

\begin{document}


\newtheorem{theorem}{Theorem}[section]
\newtheorem{definition}[theorem]{Definition}
\newtheorem{corollary}[theorem]{Corollary}
\newtheorem{lemma}[theorem]{Lemma}
\newtheorem{proposition}[theorem]{Proposition}
\newtheorem{assumption}[theorem]{Assumption}
\newtheorem{example}[theorem]{Example}
\newtheorem{remark}[theorem]{Remark}
\newcommand{\qed}{\hfill\rule{2mm}{3mm}\vspace{4mm}}

\numberwithin{equation}{section}

\def\ep{{\varepsilon}}\def\si{{\sigma}}\def\la{{\lambda}}
\def\al{{\alpha}}\def\be{{\beta}}\def\de{{\delta}}
\def\th{{\theta}}\def\ze{{\zeta}}

\def\({\left(}\def\){\right)}
\def\EE{\mathbbm{E}}
\def\PP{\mathbbm{P}}
\def\RR{\mathbbm{R}}
\def\NN{\mathbbm{N}}
\def\QQ{\mathbbm{Q}}
\def\ZZ{\mathbbm{Z}}

\def\us{{\underline{s}}}
\def\ut{{\underline{t}}}
\def\uu{{\underline{u}}}

\maketitle
\begin{abstract}
This paper studies a nonlinear filtering problem over an infinite time interval. The signal to be estimated is driven by a stochastic partial differential equation involves unknown parameters. Based on discrete observation, strongly consistent estimators of the parameters are derived at first. With the optimal filter given by Bayes formula, the uniqueness of invariant measure for the signal-filter pair has been verified. The paper then establishes approximation to the optimal filter, showing that the pathwise average distance, per unit time, of the computed approximating filter from the optimal filter converges to zero in probability. Simulation results are presented at last.

\vspace{.5cm}\noindent


\end{abstract}

\section{Introduction}
This work is concerned with approximations to a nonlinear filtering problem over an infinite time interval, which is motivated by environment pollution detection. On a river in a range of $[0,1]$, some undesired chemical are released from an unknown location $\beta$ modeled by a Poisson process with unknown intensity $\lambda$.
Let $U(t, x)$ denote the concentration of undesired chemical at time $t$ and location $x$. Then by Kallianpur and Xiong \cite{Kallianpur and Xiong 1995}, $U(t,x)$ is governed by the following stochastic partial differential equation (SPDE):
\begin{equation}\label{signal}
\partial_t U(t,x)= a \partial^2_x U(t, x)+b\partial_x U(t,x)
-\alpha U(t,x)+N_{\lambda t}\delta_{\beta}(x)
\end{equation}
given the initial $U(0, \cdot)$ and Neumann boundary conditions:
$$\partial_x U(t,0)=0, \quad \partial_x U(t,1)=0, \qquad t\ge 0, \; x\in[0,1].$$
Here $a$ is the dispersion coefficient, $b$ is the velocity of the river flow,
$\alpha$ is the dilution rate of the chemical,
$N_t$ is a standard Poisson process independent of $U(0, \cdot)$
and $\delta_{\beta}$ is the Dirac measure at $\beta$. 

The concentration of undesired chemical can rarely be observe directly. Instead, an observation station is set up and the chemical concentration around station is measured subject to random error. Herein, we consider the following observation model:
\begin{equation}\label{obs}
Y_i=h(\langle U_{t_i},\phi_0\rangle)+\varepsilon_i, \qquad i=1, 2, \cdots
\end{equation}
where $\langle U_{t_i},\phi_0\rangle=\int^1_0 \phi_0 (x)U(t_i,x)dx$.
 $h$ is a bounded continuous function on $\RR$.
 $\phi_0$ can be either bounded continuous function or Dirac delta function on $[0,1]$.
$\varepsilon'_i$s are independent with the standard normal
distribution and independent of $N_t$, $t_i=i\Delta$ with $\Delta >0$.  

Let $\mathcal{X}$ be the space of finite measures on $[0,1]$. We have the distribution valued process $U_{t}\equiv U(t,\cdot)\in\mathcal{X}$. The observation filtration ${\cal F}^Y_j=\sigma \{Y_{i}; i=1,\cdots,j\}$ contains all the available information about $\beta,\lambda,U_t$. The primary aim of our filtering problem is to get an estimate of $\beta,\lambda$ based on ${\cal F}^Y_j$, and approximate the optimal filter
\begin{equation}\label{defopt}
\langle\Pi_j,\phi\rangle=\EE\{\phi(U_{t_j})|{\mathcal{F}}^Y_j\},\quad \forall\phi\in C_b({\mathcal{X}}).
\end{equation}
It can be easily to show that the condition expectation $\EE\{\phi(U_{t_j})|{\mathcal{F}}^Y_j\}$ is the best estimate of $\phi(U_{t_j})$ based on ${\cal F}^Y_j$. 

Nonlinear filtering is widely used in engineering, natural science, machine learning and even mathematical finance; see for Van Leeuwen \cite{Van Leeuwen 2010}, Bishop \cite{Bishop 2006}, Brigo and Hanzon \cite{Brigo and Hanzon 1998}.
For the approximations to nonlinear filtering problems over the past few decades, some approaches have been developed. The pioneering work is Markov chain approximation raised by Kushner \cite{Kushner 1977} in 1977. The idea is to find a Markov chain which suitably approximated the actual signal process, then to construct the approximating filter for the Markov chain and actual observations. 
Another algorithm to solve the nonlinear filtering problems is the particle system introduced by Gordon et al. \cite{Gordon 1993} in 1993. The solution is represented by a system of weighted particles whose positions and weights are governed by SDEs that can be solved numerically; see Chapter 8 of Xiong \cite{Xiong 2008} for details. 
There are also some other numerical schemes coming out meanwhile, such as projection filter and assumed density approach; see for Brigo et al. \cite{Brigo 1999} and Kushner \cite{Kushner 2008}.
To approximate the signal process and optimal filter for model $(\ref{signal})-(\ref{obs})$, the current paper will construct Markov processes with the idea of Markov chain approximation and the obtained estimators of $\beta,\lambda$. Note that the problem is of interest over a very long time interval. Thus when talking about the error of approximated filter from optimal filter, we would use the pathwise average per unit time rather than the mean value; following the remarkable works done by Budhiraja and Kushner \cite{Budhiraja and Kushner 1999, 
Budhiraja and Kushner 2001}. Kouritzin et al \cite{Kouritzin 2004} also studied a nonlinear filtering problem over a long time interval. Their signal process was a reflecting diffusion in a random environment. In this article, the signal to be estimated is driven by an SPDE involved with jump process and unknown parameters, leading to a different way to construct its approximating filter process. 

It should be highlighted that this article will prove in a systematic way that the signal-filter pair has a unique probability invariant measure. The critical importance of uniqueness of the joint process was first raised by Kushner and Huang \cite{Kushner and Huang 1986}. 
In this paper, Theorem \ref{invariant prob}, \ref{uniqueness} will present a rigorous proof that the signal process has a unique invariant measure and the filter forgets its initial condition, yielding uniqueness of invariant measure of optimal filter by Theorem 7.1 in \cite{Budhiraja and Kushner 1999}.

The article is organized as follows. In the next section, we are to look for the estimators of $\lambda$ and $\beta$ based on the data set $\{ Y_i;i=1,\cdots,n\}$. Section 3 will construct an approximating filter with the obtained estimators and discuss the pathwise average errors over an infinite time interval. Simulation for the estimators as well as approximating filter are given in Section 4.

\section{Statistical estimation}
This section is focus on the moment estimate of parameters $\be,\la$. To make the calculation, discussion about the $r$-th sample moment is required.
Let $\mathcal{P}(\mathcal{X})$ be the collection of all probability measures on $\mathcal{X}$. Denote the law of $U_t$ by $\mathcal{L}(U_t)$. 
Here is a theorem talking about the convergence of  $\mathcal{L}(U_t)$, with which the limits of sample moments will be derived.    

\begin{theorem}\label{invariant prob}As $t\to\infty$, the sequence of measures $\mathcal{L}(U_t)$ converges in $\mathcal{P}(\mathcal{X})$ to a probability measure $F$ which is the unique invariant probability measure of the process $U_t$. Especially,  for any bounded continuous function $f$, we have
\begin{equation}\label{SLLN}
\frac 1n\sum^n_{i=1}f(Y_i)
\overset{a.s.}{\longrightarrow} \int_{\cal X}\int_{\RR} f(h(\langle \nu, \phi_0\rangle)+x)\frac{1}{\sqrt{2\pi}}e^{-\frac{x^2}{2}}F(d\nu)dx,\quad \mbox{as }n\to\infty.
\end{equation}
\end{theorem}
\noindent{\bf Proof.} From Theorem 7.2.1 of Kallianpur and Xiong \cite{Kallianpur and Xiong 1995},
equation ($\ref{signal}$) has a unique solution given by
\begin{equation}
U(t, x)=
e^{-\alpha t}\int^1_0 P_t(y, x)U_0(dy)
+\int^t_0 e^{-\alpha (t-s)}P_{t-s}(\beta, x)dN_{\lambda s}
\end{equation}
where $P_t(y,x)$ is the Green's function of the operator $L=a\partial^2_xU+b\partial_xU$ with Neumann boundary conditions given by
\[
P_t(y, x)=\sum^\infty_{j=0}e^{-\sigma_j t}\psi_j(y)\psi_j(x)
\]
and
\[
\sigma_0=0, \quad \sigma_j=a\(c^2+(j\pi)^2\);
\]
\[
\psi_0(x)=\sqrt{\frac{2c}{1-e^{-2c}}}, \quad
\psi_j(x)=\sqrt{2}e^{cx}\sin(j\pi x+k_j);
\]
\[
k_j=\tan^{-1}\big(-\frac{j\pi}{c}\big), \quad
c=-\frac{b}{2a}, \quad j=1,2,\cdots
\]

Now, we prove the tightness of $\{\mathcal{L}(U_t)\}$ in $\mathcal{P}(\mathcal{X})$. To this end, it is enough to prove that for any bounded continuous function $f$ on $[0,1]$, the family $\{\mathcal{L}(\left<U_t,f\right>)\}$ is tight in $\mathcal{P}(\RR)$. Without loss of generality, take $0\le f\le K$, and write $T_tf(y)=\int^1_0 P_t(y,x)f(x)dx$. Then
\begin{eqnarray*}
\EE\left<U_t,f\right>&=&e^{-\al t}\int^1_0 T_t f(y) U_0(dy)+\EE \int^t_0e^{-\al(t-s)}T_{t-s} f(\beta) dN_{\lambda s}\\
&\le& Ke^{-\al t}\int^1_0U_0(dy)+\int^t_0e^{-\al(t-s)}K\la ds\\
&\le&K\int^1_0U_0(dy)+K,
\end{eqnarray*}
which indicates the tightness of $\{\mathcal{L}(\left<U_t,f\right>)\}$. By Prokhorov theorem, tightness implies relative compactness, meaning that each subsequence of $\{\mathcal{L}(U_t)\}$ contains a further subsequence converging weakly. In this way, if the characteristic functions to all subsequences of $\{\left<U_t,f\right>\}$ converge to a unique function, then one can conclude that all subsequences of $\{\mathcal{L}(\left<U_t,f\right>)\}$ converges to the same distribution. In fact,
\begin{eqnarray*}
\EE\exp\(i\left<U_{t},f\right>\)&=&\EE\exp\(ie^{-\al t}\int^1_0T_tf(y)U_0(dy)
+i\int^t_0e^{-\al(t-s)}T_{t-s}f(\be)dN_{\la s}\)\\
&=&\exp\(ie^{-\al t}\int^1_0T_tf(y)U_0(dy)\)\\
&&\qquad \times\exp\(\la\int^t_0\(\exp\(ie^{-\al(t-s)}T_{t-s}f(\be)\)-1\)ds\)\\
&=&\exp\(ie^{-\al t}\int^1_0T_tf(y)U_0(dy)\)
 \exp\(\la\int^t_0\(e^{ie^{-\al s}T_{s}f(\be)}-1\)ds\)\\
&\to&\exp\(\la\int^\infty_0\(e^{ie^{-\al s}T_{s}f(\be)}-1\)ds\),\mbox{ as }t\to\infty.
\end{eqnarray*}
This implies the uniqueness of the limited measure, that is, $\mathcal{L}(U_t)\to F$ as $t\to\infty$, and hence the transformation of the Markov chain $U_t$ is uniquely ergodic. In consequence, formula $(\ref{SLLN})$ follows by the Birkhoff Ergodic Theorem.
\qed

Denote $\theta_0\equiv(\beta_0,\lambda_0)$ as the true values of $\theta\equiv(\beta,\lambda)$. Suppose the parameter $\theta_0$ taking values in a bounded set
$$\Theta=\{(\beta,\lambda); 0<\beta < 1, 0<\lambda\le M\}.$$ 
Let
$$ g_1(\theta)=\int_{\cal X} h(\langle \nu, \phi_0\rangle)F(d\nu),
\qquad g_2(\theta)=\int_{\cal X} h^2(\langle \nu, \phi_0\rangle)F(d\nu),$$
and the Jacobian matrix 
$$J(\theta)=\left[\begin{array}{cc}
\partial g_1(\theta)/\partial \beta & \partial g_1(\theta)/\partial \lambda\\
 \partial g_2(\theta)/\partial \beta &\partial g_2(\theta)/\partial \lambda
\end{array}\right].$$

\begin{theorem}\label{consistency of theta}
Let $m_r(n)=n^{-1}\sum^n_{i=1}Y^r_i$, $r=1,2$.
Suppose $g_r(\theta)$ posses continuous partial derivative on $\Theta$ and 
the Jacobian $\det(J(\theta))$ be different from zero everywhere on $\Theta$.
If the following equations
$$\begin{array}{cc}
m_1(n)=g_1(\theta)\\
m_2(n)=g_2(\theta)+1
\end{array}$$
have a unique solution, denoted by
 $\hat\theta_n=(\hat\beta_n, \hat\lambda_n)$,
with probability $1$ as $n\rightarrow\infty$,
then this solution is a strongly consistent estimator for $(\beta_0, \lambda_0)$.
\end{theorem}
{\bf Proof.} Let $S$ be the image of the set $\Theta$ under the continuous mapping $g=(g_1,g_2+1)$. With probability 1 the inverse mapping $g^{-1}$ is locally one-to-one and is a continuous mapping from $S$ into $\Theta$. 
Denote $m(n)=(m_1(n),m_2(n))$. From Theorem $\ref{invariant prob}$, we have $$m(n)\overset{a.s.}{\to}g(\theta_0),\quad \mbox{as }n\to\infty.$$ 
Hence for large $n$ with probability 1, the point $m(n)\in S$; in this case, equations $m(n)=g(\theta)$ are solvable and their solution is necessarily of the form $\hat\theta_n=g^{-1}(m(n))\in\Theta$. As a result, thanks to the continuity of $g^{-1}$, we have 
$$\hat\theta_n=g^{-1}(m(n))\overset{a.s.}{\to}\theta_0,\quad \mbox{as }n\to\infty,$$ 
and finish the proof.
\qed

With the obtained estimators, now we are able to construct the Markov processes to efficiently approximate the optimal filter $\Pi_\cdot$ in (\ref{defopt}). 

\section{Stochastic filtering}

Let ${\cal F}^N_j=\sigma \{N_u; 0\leq u\leq t_j,t_j=j\Delta\}$
and ${\cal F}_j={\cal F}^N_{j}\vee {\cal F}^Y_j$. We take simplified notations $Y_{0,j}= \{Y_{i}; i=1,\cdots,j\}$ and $U_{0,j}=\{U_{t_i};i=1,\cdots,j\}$.
For each $j\in \NN$, let
$$
R(U_{0,j},Y_{0,j})=\exp\left\{\sum^j_{i=1}[Y_{i}h(\langle U_{t_i},\phi_0\rangle)
-\frac 12h^2(\langle U_{t_i},\phi_0\rangle)]\right\}.$$
Then, one can define a probability measure $\QQ$ by
$$\frac{d\QQ}{d\PP}\left |_{{\cal F}_j}=R^{-1}(U_{0,j},Y_{0,j})\right..$$
Let $\EE_{\QQ}$ denote the expectation with respect to $\QQ$. Unless further notice, $\EE$ denote the expectation with respect to the original probability measure $\PP$.
In the following theorem, we have an expression for the optimal filter $\Pi_\cdot$ in terms of $\QQ$.
\begin{theorem}\label{indep}
\begin{itemize}
\item[$(\romannumeral1)$] Under $\QQ$, $\{Y_1, \cdots, Y_j\}$ are mutually independent
standard normal random variables and are independent of
$(N_s)_{0\le s\le t_j}$.
\item[$(\romannumeral2)$] Under $\QQ$, $(N_s)_{0\le s\le t_n}$ has the same distribution
as under $\PP$.
\item[$(\romannumeral3)$] (The Bayes formula) For any bounded and continuous
function $\phi$,
\begin{eqnarray}\label{bayes}
\EE\big\{\phi(U_{t_j})|{\mathcal{F}}^Y_{j}\big\}
&=&\frac {\EE_{\QQ}\big\{\phi(U_{t_j})
\frac {d\PP}{d\QQ}|{\mathcal{F}}^Y_{j}\big\}}
{\EE_{\QQ}\big\{\frac {d\PP}{d\QQ}|{\mathcal{F}}^Y_{j}\big\}} \nonumber\\
&=&\frac {\EE_{\QQ}\big\{\phi(U_{t_j})R(U_{0,j},Y_{0,j})|{\mathcal{F}}^Y_{j}\big\}}
{\EE_{\QQ}\big\{R(U_{0,j},Y_{0,j})|{\mathcal{F}}^Y_{j}\big\}}.
\end{eqnarray}
\end{itemize}
\end{theorem}
The proof is done by the similar arguments used in Theorem 3.1, 3.2 of Mandal and Mandrekar \cite{Mandal and Mandrekar 2000}. 

Let $\tilde U_0$ be a measure-valued random variable
with the same distribution as $U_0$ but to be independent of
$U_0$ and all other processes. Let $\tilde N$ be a standard Poisson process
independent of all other processes. Define
\begin{equation}
\tilde U(t,x)=e^{-\alpha t}\int^1_0P_t(y,x)\tilde U_0(dy)
+\int^t_0 e^{-\alpha (t-s)}
P_{t-s}(\beta_0, x)d\tilde N_{\lambda_0 s}.
\end{equation}
Thanks to Theorem $\ref{indep}$, the condition expectation $(\ref{bayes})$ can be rewritten as
\begin{equation*}
\langle\Pi_j,\phi\rangle
=\frac {\EE_{\{\Pi_0,Y_{0,j}\}}\big\{\phi(\tilde U_{t_j})
R(\tilde U_{0,j},Y_{0,j})\big\}}
{\EE_{\{\Pi_0,Y_{0,j}\}}\big\{R(\tilde U_{0,j},Y_{0,j})\big\}},
\end{equation*}
where $\EE_{\{\Pi,Y\}}G(\tilde U_t, Y)$ denotes the conditional expectation of
a function $G$ of $\tilde U_t$ and $Y$, given the data $Y$ and the initial distribution of $\tilde U_t$ is $\Pi$.
In the end, owing to the Markov property of $U_t$,
the optimal filter satisfies the semigroup relation
\begin{equation}\label{opt}
\langle\Pi_j,\phi\rangle
=\frac {\EE_{\{\Pi_{j-1},Y_{j}\}}\big\{\phi(\tilde U_{t_1})
R(\tilde U_{t_1},Y_{j})\big\}}
{\EE_{\{\Pi_{j-1},Y_{j}\}}\big\{R(\tilde U_{t_1},Y_{j})\big\}},
\end{equation}
on which we impose a definition $R(\tilde U_{t_1},Y_{j})=\exp\big\{Y_{j}h(\langle \tilde U_{t_1},\phi_0\rangle)-\frac{1}{2}h^2(\langle \tilde U_{t_1},\phi_0\rangle)\big\}$. In $(\ref{opt})$, $\Pi_{j-1}$ is the distribution of $\tilde U_0$.

The following is a critical important theorem of this paper.
\begin{theorem}\label{uniqueness}
The signal-filter pair $\{U_{t_j}, \Pi_j\}$ has a unique invariant probability measure.
\end{theorem}

\noindent{\bf Proof.}  Note that $U_{t_j}$ is a hidden Markov process. We first verify that
\begin{equation}\label{unique}
\EE\(\|\PP(U_{t_j}\in\cdot|U_0)-\mu_0\|_{TV}\)\to 0,\qquad\mbox{ as }j\to\infty,
\end{equation}
with invariant measure $\mu_0\equiv\PP(V_{\infty}\in\cdot)$, where $\|\cdot\|_{TV}$ is the total variation norm and 
\[V_j=\int^{t_j}_0e^{-\al (t_j-s)}P_{t_j-s}(\be,\cdot)dN_{\la s}.\]
Note that for each $j$, $\{N_{\la s},\,0\le s\le t\}$ and $\{N_{\la t_j}-N_{\la(t_j-s)}\}$ are equal in distribution and hence,
\begin{eqnarray*}
V_j&\stackrel{\mathcal{D}}{=}&\int^{t_j}_0e^{-\al (t_j-s)}P_{t_j-s}(\be,\cdot)d\(N_{\la t_j}-N_{\la(t_j-s)}\)\\
&=&\int^{t_j}_0e^{-\al s}P_{s}(\be,\cdot)dN_{\la s}\\
&\to&V_\infty\equiv\int^\infty_0e^{-\al s}P_{s}(\be,\cdot)dN_{\la s}.
\end{eqnarray*}
Note that $U_{t_j}-V_j\to 0$. Let $(U_{t_j})$ be extended to $j\in\ZZ$ as a stationary process and let $\PP^x$ be a version of regular conditional probability $\PP^{U_0}=\PP(\cdot|U_0)$. Since the limit of $U_{t_j}$ does not depend on the initial, it is easy to verify that for any $A\in\cap_{j\le 0}\mathcal{F}^U_{-\infty,j}$ and $x,y\in\mathcal{X}$, $\PP^x(A)=\PP^y(A)\in\{0,1\}$. It then follows from the remark after Conjecture 1.7 of van Handal \cite{Handel 2012} or Theorem 2.3 of van Handel \cite{Handel 2009} that (\ref{unique}) holds.

Further, the non-degenerate property (cf. Definition 1.5 of \cite{Handel 2012}) of the observations follows from
(\ref{obs}) directly. Hence, by Theorem 1.6 of  \cite{Handel 2012} that the exchangeability (1.1) of that paper holds.

One can see that the limitation of $U_t$, as $t$ goes to $\infty$,
is independent of the initial value $U_0$, and so the filter
$\Pi_j$ forgets its initial condition. Thus, the signal-filter pair
$\{U_{t_j}, \Pi_j\}$ has a unique probability invariant measure
by Theorem 7.1 of Budhiraja and Kusher \cite{Budhiraja and Kushner 1999}.
\qed

Here we proceed to the establishments of approximating signal and the computed approximating filter. Set
\begin{equation}
\tilde U^n(t, x)=e^{-\alpha t}\int^1_0P_t(y,x)\tilde U_0(dy)
+\int^t_0 e^{-\alpha (t-s)}
P_{t-s}(\hat\beta_n, x)d\tilde N_{\hat\lambda_n s},
\end{equation}
and $\tilde U^n_t\equiv \tilde U^n(t, \cdot)$. Define the approximating filter $\Pi^n$ recursively by
\begin{equation}\label{app fil}
\langle\Pi^n_j,\phi\rangle
=\frac {\EE_{\{\Pi^n_{j-1},Y_{j}\}}\big\{\phi(\tilde U^n_{t_1})
R(\tilde U^n_{t_1},Y_{j})\big\}}
{\EE_{\{\Pi^n_{j-1},Y_j\}}\big\{R(\tilde U^n_{t_1},Y_{j})\big\}},\quad n\le j.
\end{equation}
The next theorem shows the convergence property of approximating signal $\tilde U^n_{t_1}$.
\begin{theorem}\label{convergence of U}
For any bounded continuous function $\phi$,
$$\EE\big|\langle\tilde U^n_{t_1}, \phi\rangle-
\langle\tilde U_{t_1}, \phi\rangle\big|^2 \to 0,
\qquad as \ n \rightarrow \infty,$$
where $\Pi^n_{j-1}$, $j=1,2,\cdots$ is the distribution of $\tilde U_0$.
\end{theorem}
\noindent{\bf Proof.}
Note that for any bounded continuous function $\phi$,
\begin{align*}
&\EE\big|\langle\tilde{U}^n_{t_1},\phi\rangle-\langle\tilde{U}_{t_1},\phi\rangle\big|^2 \\
=&\EE\Big|\int^{t_1}_0 e^{-\alpha (t_1-s)}
\int_0^1 \phi(x)P_{t_1-s}(\hat\beta_n, x)dxd\tilde N_{\hat\lambda_n s}\\
& ~~~~~~-\int^{t_1}_0 e^{-\alpha (t_1-s)}
\int_0^1 \phi(x)P_{t_1-s}(\beta_0, x)dxd\tilde N_{\lambda_0 s}\Big|^2\\
\le&2\EE \Big|\int^{t_1}_0 e^{-\alpha (t_1-s)}
\int_0^1  \phi(x)\Big[P_{t_1-s}(\hat\beta_n, x)dx
-P_{t_1-s}(\beta_0, x)\Big]dxd \tilde N_{\hat\lambda_n s}\Big|^2\\
&~~~+2\EE\Big|\int^{t_1}_0 e^{-\alpha (t_1-s)}
\int_0^1 \phi(x)P_{t_1-s}(\beta_0, x)dx
d\Big[\tilde N_{\hat\lambda_n s}-\tilde N_{\lambda_0 s}\Big]\Big|^2\\
\le&K\EE\hat\lambda_n(1+t_1 \hat\lambda_n)\int^{t_1}_0 e^{-2\alpha (t_1-s)}
\int_0^1  \Big|P_{t_1-s}(\hat\beta_n, x)dx
-P_{t_1-s}(\beta_0, x)\Big|^2dxds\\    
&~~~+2\int^{t_1}_0 e^{-2\alpha (t_1-s)}
\Big(\int_0^1 \phi(x)P_{t_1-s}(\beta_0, x)dx\Big)^2
\EE\big|\hat\lambda_n-\lambda_0\big|\big(1+t_1\big|\hat\lambda_n-\lambda_0\big|\big)ds.
\end{align*}
Obviously, the second term on the RHS goes to zero with Theorem $\ref{consistency of theta}$. For the double integral,
\begin{align*}
&\int^{t_1}_0\int_0^1  \Big|P_{t_1-s}(\hat\beta_n, x)dx
-P_{t_1-s}(\beta_0, x)\Big|^2dxds\\
\le&\int_0^1 2\sum_{j=0}^\infty \int^{t_1}_0 e^{-2\sigma_j(t_1-s)}ds\psi^2_j(x)
\Big|e^{c\hat\beta_n}\sin(j\pi\hat\beta_n+k_j)
-e^{c\beta_0}\sin(j\pi\beta_0+k_j)\Big|^2 dx\\
\le& 8e^{2c}\sum_{j=0}^\infty \frac{1}{2\sigma_j}(1-e^{-2\sigma_j t_1})
\Big(\big|e^{c\hat\beta_n}-e^{c\beta_0}\big|^2\\
&~~~~~~~~~~~
+2e^{2c\beta_0}\Big|\sin\frac{j\pi(\hat\beta_n-\beta_0)}{2}
\cos\frac{j\pi(\hat\beta_n+\beta_0)+2k_j}{2}\Big|^2\Big)\\
\le& K\sum_{j=0}^\infty \frac{1}{\sigma_j}
\Big(\big|e^{c\hat\beta_n}-e^{c\beta_0}\big|^2
+\Big|\sin\frac{j\pi(\hat\beta_n-\beta_0)}{2}\Big|^2\Big)\\
\le& \sum_{j=0}^\infty \frac{K}{a(c^2+(j\pi)^2)}
\Big(\big(e^{c\delta^*}(\hat\beta_n-\beta_0)\big)^2
+\big(j\pi|\hat\beta_n-\beta_0|\big)^{\frac{1}{2}}\Big)\\
\le& \sum_{j=0}^\infty \frac{K}{a(c^2+(j\pi)^2)}
\Big((\hat\beta_n-\beta_0)^2
+\big(j\pi|\hat\beta_n-\beta_0|\big)^{\frac{1}{2}}\Big),
\end{align*}
where $\delta^*=\delta^*(\beta_0, \hat\beta_n)$ falls between $\beta_0$ and $\hat\beta_n$. With Theorem $\ref{consistency of theta}$, we obtain that the double integral converges to $0$ in $L^1$-norm.
\qed

We now need to investigate the asymptotic behavior of $\{U_{t_j},\Pi^n_j\}$ and prove the main result. Recall that in Theorem \ref{uniqueness}, $\{U_{t_j},\Pi_j\}$ is shown to has a unique invariant measure, denoted by $\rho$. Let $f$ be any bounded continuous real-valued function. In what follows, we will show that
\begin{equation}\label{PIn to rho}
\frac1k\sum^k_{j=1}f(U_{t_j},\Pi^n_j)
\rightarrow\int_{\mathcal{X}\times\mathcal{P}(\mathcal{X})} f(\nu,\alpha) \rho(d\nu,d\alpha)
\end{equation}
in probability as $n,k$ go to $\infty$.

Let's get start on compactifying the finite measure space $\mathcal{X}$ with the idea of one-point compactification in real space. Set
$$
\bar{\cal X}={\cal X}\bigcup\{\mu; \mu \mbox{ is any positive
infinite measure on } [0,1]\}.
$$
Then $\bar{\mathcal{X}}$ includes each compact set $B\subseteq\mathcal{X}$, and the “infinite set” of the form $(\mathcal{X}\setminus B)\cup\{\mu; \mu([0,1])>k,k\ge0\}$. Thus any sequence of measures on $\bar{\mathcal{X}}$ is tight.

For each $j$ and $n$, define
$$ B_j=\sum^j_{i=1}\varepsilon_i,\quad 
\Psi^n(j,\cdot)=(U_{t_{j+\cdot}}, \Pi^n_{j+\cdot}, Y_{j+\cdot},
B_{j+\cdot}-B_j).$$
Let $I_C(\Psi^n(j,\cdot))$ denote the indicator function of the event
$\Psi^n(j,\cdot)\in C$ and $C$ be any measurable set in the product path space of $\Psi^n(j,\cdot)$.
Here we define the occupation measures:
\begin{equation}
Q^{n,k}(C)=\frac 1 k\sum^k_{j=1}I_C\big(\Psi^n(j,\cdot)\big)
\end{equation}
that take value in the space of probability measures on the product path space
$$\mathcal{P}\big(D(\bar{\mathcal{X}};\NN)\times D(\mathcal{P}(\bar{\mathcal{X}});\NN)\times D(\RR;\NN)\times D(\RR;\NN)\big).$$
$D(E;\NN)$ is the path space of $E$-valued functions. 
To verify formula (\ref{PIn to rho}), in the sequel, we will discuss the tightness and weak convergence of the measure-valued random variables $Q^{n,k}(\cdot)$.

\begin{theorem}\label{tight} 
$\{U_{t_{j+\cdot}}; j\ge 0\}$ is tight
in $D(\bar{\mathcal{X}},\NN)$;
$\{\Pi^n_{j+\cdot}; n>0, j\ge 0\}$ is tight in
$D({\mathcal{P}}(\bar{\cal X}),\NN)$;
$\{Y_{j+\cdot}; j\ge 0\}$ is tight in $D(\RR,\NN)$;
and $\{Q^{n,k}; n>0, k>0\} $ is tight on
$\mathcal{P}\big(D(\bar{\mathcal{X}};\NN)\times D(\mathcal{P}(\bar{\mathcal{X}});\NN)\times D(\RR;\NN)\times D(\RR;\NN)\big)$.
\end{theorem}
{\bf Proof.}

$(\romannumeral1)$ $\{U_{t_{j+\cdot}}; j\ge 0\}$ is tight since $\bar{\mathcal{X}}$ is compact.

$(\romannumeral2)$ Obviously,  $\{\Pi^n_{j+\cdot}; n>0, j\ge 0\}$ is tight
in $D({\mathcal{P}}(\bar{\cal X}),\NN)$ since this space is compact.

$(\romannumeral3)$ $\{Y_{j+\cdot}; j\ge 0\}$ is tight from the boundedness of $h$
and stationary of $\varepsilon_j$.

$(\romannumeral4)$ The sequence $\{\EE[Q^{n,k}(\cdot)]\}$ is tight from the tightness
of  $\{U_{t_{j+\cdot}}; j\ge 0\}$, $\{\Pi^n_{j+\cdot}; n>0, j\ge 0\}$,
 $\{Y_{j+\cdot}; j\ge 0\}$ and $\{B_{j+\cdot}-B_j;j\ge0\}$.
Thus, the sequence $\{Q^{n,k}(\cdot)\}$ is tight
by Theorem 1.6.1 of Kushner \cite{Kushner 1990}.
\qed

Thanks to the tightness property, we understand that each subsequence of $\{Q^{n,k}(\cdot)\}$ converges weakly. Before we identify the limits of subsequences, define the occupation measures for the signal-filter pair
\begin{equation}
Q^{n,k}_1(C)=\frac 1 k\sum^k_{j=1}I_C\big(\Psi^n_1(j,\cdot)\big),\quad
\Psi^n_1(j,\cdot)=(U_{t_{j+\cdot}}, \Pi^n_{j+\cdot}).
\end{equation}
\begin{theorem}\label{stationarity}
Let $Q(\cdot)$ denote any weak sense limit of $\{Q^{n,k}(\cdot)\}$.
Let $\omega$ be the canonical variable on the probability space
on which $Q(\cdot)$ is defined, and denote the samples by $Q^{\omega}$.
For each $\omega$, $Q^{\omega}$ induces a process
$$\Psi^\omega_{\cdot}
=(U^\omega_{t_\cdot}, \Pi^{\omega}_{\cdot},
Y^\omega_\cdot, B^\omega_\cdot).$$
For almost all $\omega$, the following hold.
\begin{itemize}
\item[$(\romannumeral1)$] $\Pi^{\omega}_{\cdot}$ is ${\mathcal{P}}({\mathcal{X}})$-valued.

\item[$(\romannumeral2)$] $(U^{\omega}_{t_\cdot}, \Pi^{\omega}_{\cdot})$ is stationary.

\item[$(\romannumeral3)$] $\{B^{\omega}_{\cdot}\}$ is the sum of mutually independent
$\mathcal{N}(0,1)$ random variables $\{\varepsilon^{\omega}_i\}$
which are independent of $U^{\omega}_{\cdot}$.

\item[$(\romannumeral4)$] $Y^{\omega}_\cdot=h(\langle U^{\omega}_{t_\cdot}, \phi_0\rangle)
+\varepsilon^{\omega}_\cdot$.

\item[$(\romannumeral5)$] $U^{\omega}_{t_\cdot}$ has the transition function of $U_{t_\cdot}$.

\item[$(\romannumeral6)$]  for each integer $j$ and $\phi\in C_b({\mathcal{X}})$,
\begin{equation}\label{<Pi,phi>}
\langle\Pi^{\omega}_j,\phi\rangle
=\frac {\EE_{\{\Pi^{\omega}_0,Y^{\omega}_{0,j}\}}\big\{\phi(\tilde U_{t_j})
R(\tilde U_{0,j},Y^{\omega}_{0,j})\big\}}
{\EE_{\{\Pi^{\omega}_0,Y^{\omega}_{0,j}\}}
\big\{R(\tilde U_{0,j},Y^{\omega}_{0,j})\big\}}.
\end{equation}
\end{itemize}
\end{theorem}

\noindent{\bf Proof.} Assume the lower case $\psi_\cdot=(u_\cdot,\pi_\cdot,y_\cdot,b_\cdot)$ is the canonical sample path in the space where $\Psi^\omega_\cdot$ (of all $\omega$) are defined.
Some of the ideas to be used in the proof come from Budhiraja and Kushner \cite{Budhiraja and Kushner 1999}.
\begin{itemize}
\item[$(\romannumeral1)$] Let $K_\delta=\{\mu\in \bar{\mathcal{X}}; \mu([0,1])\le 1/\delta\}$.
By the definition of the $\Pi^n_j$, we have that
$$
\Pi^n_{j}(K^c_{\delta^j})=
\frac {\EE_{\{\Pi^n_{j-1},Y_j\}}
\big\{I_{K^c_{\delta^j}}(\tilde U^n_{t_1})
R(\tilde U^n_{t_1},Y_{j})\big\}}
{\EE_{\{\Pi^n_{j-1},Y_j\}}\big\{R(\tilde U^n_{t_1},Y_{j})\big\}}.
$$
Let $\tilde{V}^n_1=\int^{t_1}_0e^{-\alpha(t_1-s)}
P_{t_1-s}(\hat\beta_n, \cdot)d\tilde{N}_{\hat\lambda_n s}$.
Thus, we have that, for each $n$ and $j$,
\begin{align*}
&\Pi^n_{j}(K^c_{\delta^j})\\
=&\int \EE\Big\{I_{K^c_{\delta^j}}
\Big(\int e^{-\alpha t_1}P_{t_1}(y, \cdot)\nu(dy)+\tilde V^n_{1}\Big)
R\Big(\int e^{-\alpha t_1}P_{t_1}(y,\cdot)\nu(dy)+\tilde V^n_{1},Y_j\Big)\\
&~~~~~~~~\Big|Y_{j}\Big\} \Pi^n_{j-1}(d\nu)
\Big/ 
\int \EE\Big\{R\Big(\int e^{-\alpha t_1}P_{t_1}(y,\cdot)\nu(dy)
+\tilde V^n_{1},Y_{j}\Big)|Y_j\Big\}\Pi^n_{j-1}(d\nu)\\
\le&\int \EE\Big[I_{K^c_{\delta^j}}
\Big(\int e^{-\alpha t_1}P_{t_1}(y, \cdot)\nu(dy)
+\tilde V^n_{1}\Big)
\exp\Big\{h\Big(\big<\int e^{-\alpha t_1}P_{t_1}(y,\cdot)\nu(dy)\\
&~~~~~~+\tilde V^n_{1}, \phi_0\big>\Big) Y_{j} 
-\frac{1}{2}h^2\Big(\big<\int e^{-\alpha t_1}P_{t_1}(y,\cdot)\nu(dy)
+\tilde V^n_{1}, \phi_0\big>\Big)\Big\}
\Big| Y_{j}\Big]\Pi^n_{j-1}(d\nu) \\
&~~~~~~\Big/\Big\{\int \EE\Big[\exp\Big\{h\Big(\big<\int e^{-\alpha t_1}P_{t_1}(y, \cdot)\nu(dy)
+\tilde V^n_{1}, \phi_0\big>\Big) Y_{j}\\
& ~~~~~~-\frac{1}{2}h^2\Big(\big<\int e^{-\alpha t_1}P_{t_1}(y,\cdot)\nu(dy)
+\tilde V^n_{1}, \phi_0\big>\Big)\Big\}
\Big|Y_{j}\Big]\Pi^n_{j-1}(d\nu)\Big\}\\
\le&c_1e^{c_2 |\varepsilon_{j}|}\int \EE\Big[I_{K^c_{\delta^j}}
\Big(\int e^{-\alpha t_1}P_{t_1}(y, \cdot)\nu(dy)
+\tilde V^n_{1}\Big)\Big|Y_{j}\Big]\Pi^n_{j-1}(d\nu)\\
=&c_1e^{c_2 |\varepsilon_{j}|}\Big(\int_{K_{\delta^{j-1}}}
+\int_{K^c_{\delta^{j-1}}}\Big)
\EE\Big\{I_{K^c_{\delta^j}}
\Big(\int e^{-\alpha t_1}P_{t_1}(y, \cdot)\nu(dy)
+\tilde V^n_{1}\Big)\Big|Y_{j}\Big\}\Pi^n_{j-1}(d\nu)\\
=&c_1e^{c_2 |\varepsilon_{j}|}\Big(\int_{K_{\delta^{j-1}}}
+\int_{K^c_{\delta^{j-1}}}\Big)
\PP\Big\{\Big(\int e^{-\alpha t_1}P_{t_1}(y, \cdot)\nu(dy)
+\tilde V^n_{1}\Big)\ge \frac{1}{\delta^j}\Big|Y_{j}\Big\}\Pi^n_{j-1}(d\nu)\\
\le&c_1 e^{c_2 |\varepsilon_{j}|}
\Big\{\delta^j\int_{K_{\delta^{j-1}}}
\EE\left[e^{-\alpha t_1}\nu([0,1])
+\tilde V^n_{1}\big|Y_{j}\right]
\Pi^n_{j-1}(d\nu)+\Pi^n_{j-1}(K^c_{\delta^{j-1}})\Big\}\\
\le&c_1 e^{c_2 |\varepsilon_{j}|}
[\delta(1+\hat \lambda_n)+\Pi^n_{j-1}(K^c_{\delta^{j-1}})]\\
\le&\delta(1+\hat \lambda_n)\sum^j_{k=2}\prod^j_{i=k}(c_1e^{c_2 |\varepsilon_{i}|})
+\prod^j_{i=2}(c_1e^{c_2 |\varepsilon_{i}|})\Pi^n_1(K^c_\delta)\\
\le&\delta(1+\hat \lambda_n)\sum^j_{k=2}\prod^j_{i=k}(c_2e^{c_1 |\varepsilon_{i}|})
+\delta(1+\hat \lambda_n)\prod^j_{i=1}(c_1e^{c_2 |\varepsilon_{i}|})
\int\nu[0,1]\pi_0(d\nu),
\end{align*}
since $h(\cdot)$ is bounded, where $c_1$ and $c_2$ are finite constants.
So,
$$\EE[\Pi^n_{j}(K^c_\delta)]\le c_j\delta.$$
Let $\delta\rightarrow 0 $, we have that, for almost all $\omega$,
$$\Pi^n_j\{\mu; \mu([0,1])=\infty\}=0.$$
Therefore, the proof of $(\romannumeral1)$ is completed.

\item[$(\romannumeral2)$] 
Let $H$ be the Borel set in the product $(u_\cdot,\pi_\cdot)$ path space. For arbitrary integer $c>0$, define the left shift $H_c$ by $H_c=\{(u_\cdot,\pi_\cdot); (u_{c+\cdot},\pi_{c+\cdot})\in H\}$. Then
$$
Q^{n,k}_1(H_c)=\frac1k \sum^k_{j=1} I_{H_c}\big(\Psi^n_1(j,\cdot)\big)
=\frac1k\sum^k_{j=1} I_{H}\big(\Psi^n_1(j+c,\cdot)\big),
$$
and hence 
\begin{align*}
\big|Q^{n,k}_1(H_c)-Q^{n,k}_1(H)\big|
=&\Big|\frac1k\sum^{k+c}_{j=c+1} I_{H}\big(\Psi^n_1(j,\cdot)\big)-
\frac1k\sum^k_{j=1} I_{H}\big(\Psi^n_1(j,\cdot)\big)\Big|\\
=&\Big|\frac1k \Big(\sum^{k+c}_{j=k+1}-\sum^{c}_{j=1}\Big) I_{H}\big(\Psi^n_1(j,\cdot)\big)\Big|
\le\frac{2c}{k}\rightarrow 0,
\end{align*}
as $k\to \infty$ for all $\omega,c,n,H$. In addition to the weak convergence, it yields the stationarity of their limited processes $\big(U^\omega_\cdot,\Pi^\omega_\cdot\big)$ for almost all $\omega$.

\item[$(\romannumeral3)$] Let $\phi_r(\cdot),\varphi_p(\cdot)$ be two arbitrary finite collections of bounded and continuous real-valued functions. Let $c_i$ be a finite collection of arbitrary integers.
For arbitrary integers $m$, $k$ and bounded and continuous function $q(\cdot)$,
define the function
\begin{align*}
f(\psi_\cdot)=&q\big(\langle u_{c_i},\phi_r\rangle, \langle\pi_{c_i},\varphi_p\rangle,
b_{c_i}; c_i\le m, r,p\big)\\
&\times\Big[\exp\Big\{\sum^{m+k}_{l=m}\langle u_l,\phi_l\rangle (b_{l+1}-b_l)\Big\}
-\exp\Big\{\frac{1}{2}\sum^{m+k}_{l=m}\langle u_l,\phi_l\rangle^2\Big\}\Big].
\end{align*}
We need to show that, for almost all $\omega$,
$$
\int f(\psi_\cdot)Q^{\omega}(d\psi)=0.  \eqno{\textcircled{1}}
$$
This, in turn, implies that (for almost all $\omega$) the $\varepsilon_l^{\omega}$
are mutually independent $\mathcal{N}(0,1)$ random variables,
and are independent of the $U^{\omega}_{\cdot}$ and of the ``past'' of the
$\Pi^{\omega}_{\cdot}$ process. 
To get $\textcircled{1}$, it is enough to show that as $n,k\to\infty$,
$$
\EE\Big[\int f(\psi_\cdot)Q^{n,k}(d\psi)\Big]^2
=\EE\Big[\frac1k\sum^k_{j=1}f\big(\Psi^n(j,c_i)\big)   \Big]^2
\rightarrow0. 
$$
In fact, by the definition of $\Psi^n(j,\cdot)$,
\begin{align*}
&f\big(\Psi^n(j,c_i)\big)=q\big(\langle U_{t_j+c_i},\phi_j\rangle, 
\langle\Pi^n_{j+c_i},\varphi_p\rangle,
B_{j+c_i}-B_j; c_i\le m, j,p\big)\\
&~~~~~~~~~\times\Big[\exp\Big\{\sum^{m+k}_{l=m}\langle U_{t_j+l},\phi_l\rangle 
(B_{j+l+1}-B_{j+l})\Big\}
-\exp\Big\{\frac{1}{2}\sum^{m+k}_{l=m}\langle U_{t_j+l},\phi_l\rangle^2\Big\}\Big].
\end{align*}
Using the fact that $B_{j+l+1}-B_{j+l}=\ep_{j+l}\sim\mathcal{N}(0,1)$ are mutually independent among all the $j,l$, so one can see that the mean square converges to 0.

\item[$(\romannumeral4)$] 
Write $|z|^2_1=\min\{|z|^2,1\}$. To identify $Y^\omega_\cdot$, it is necessary to prove for each integer $s$,
$$
\EE\int Q^\omega(d\psi)\big|y_s-h\big(\langle u_s,\phi_0\rangle\big)-(b_{s+1}-b_s)\big|^2_1=0.
\eqno{\textcircled{2}}
$$
which implies $(\romannumeral4)$ for almost all $\omega$. This can be shown by using the definition of occupation measures and the weak convergence. Note that
\begin{eqnarray*}
&&\EE\int Q^{n,k}(d\psi)\big|y_s-h\big(\langle u_s,\phi_0\rangle\big)-(b_{s+1}-b_s)\big|^2_1\\
&=&\EE\frac1k\sum^k_{j=1}\Big|Y_{j+s}-h\big(\langle U_{t_j+s},\phi_0\rangle\big)
-(B_{j+s+1}-B_{j+s})\Big|^2_1=0.
\end{eqnarray*}
Then, since the integrand in $\textcircled{2}$ is a bounded and continuous function of $\psi_\cdot$ and $Q^{n,k}$ weakly converges to $Q$, equality $\textcircled{2}$ follows.

\item[$(\romannumeral5)$] Adopt the notations in $(\romannumeral3)$. Using the similar arguments in \cite{Budhiraja and Kushner 1999},
we can show that, for any real number $x$,
$$
\int Q^{\omega}(d\psi)q(\langle u_{c_i}, \phi_r\rangle; c_i\le m, r)
\left[e^{ix\langle u_{m+1}, \phi\rangle}-\int e^{ix\langle \mu, \phi\rangle}
\PP(d\mu,1|u_m)\right]=0,
$$
where $\PP(d\mu,1|x)$ is the one step transition function of $U_{\cdot}$.
This implies that $U^{\omega}_{\cdot}$ is Markov with the transition
function of $U_{\cdot}$.

\item[$(\romannumeral6)$] For arbitrary canonical sample path $\psi_\cdot$ and integer $m$, define a function 
$$
A(\psi_m)=\langle\pi_m,\phi\rangle
-\frac {\EE_{\{\pi_{m-1},y_m\}}\{\phi(\tilde U_{t_1})
R(\tilde U_{t_1},y_m)\}}
{\EE_{\{\pi_{m-1},y_m\}}\{R(\tilde U_{t_1},y_m)\}}.
$$
We are to prove that, for almost all $\omega$ and all $m$,
$A(\Psi^{\omega}_m)=0$,
with probability $1$ which implies $(\ref{<Pi,phi>})$. This can be done by showing that
$$
0=\EE\int Q^{\omega}(d\psi)\big|A(\psi_m)\big|^2_1.  \eqno{\textcircled{3}}
$$
By the weak convergence, the prelimit form of the right side of $\textcircled{3}$
is $$\EE\int Q^{n,k}(d\psi)\big|A(\psi_m)\big|^2_1
=\frac{1}{k}\sum^k_{j=1}\EE\big|A(\Psi^n(j,m))\big|^2_1.$$
Thus, to obtain $\textcircled{3}$ it suffices to show that
$\EE\big|A(\Psi^n(j,\cdot)\big|^2_1\rightarrow 0$ uniformly in $j$ as $n\rightarrow \infty$.
Actually, by Definition \ref{app fil}, Theorem $\ref{convergence of U}$ and the continuity of $\phi,\phi_0,h$, we have
\begin{eqnarray*}
\EE\big|A(\Psi^n(j,\cdot))\big|^2_1
&=&\EE\Big|\langle\Pi^n_j,\phi\rangle
 -\frac {\EE_{\{\Pi^n_{j-1},Y_j\}}\{\phi(\tilde U_{t_1})
R(\tilde U_{t_1},Y_j)\}}
{\EE_{\{\Pi^n_{j-1},Y_j\}}\{R(\tilde U_{t_1},Y_j)\}}\Big|^2_1\\
&=&\EE\Big|\frac {\EE_{\{\Pi^n_{j-1},Y_j\}}\{\phi(\tilde U^n_{t_1})
R(\tilde U^n_{t_1},Y_j)\}}
{\EE_{\{\Pi^n_{j-1},Y_j\}}\{R(\tilde U^n_{t_1},Y_j)\}}\\
&&~~~~~~~~~-\frac {\EE_{\{\Pi^n_{j-1},Y_j\}}\{\phi(\tilde U_{t_1})
R(\tilde U_{t_1},Y_j)\}}
{\EE_{\{\Pi^n_{j-1},Y_j\}}\{R(\tilde U_{t_1},Y_j)\}}\Big|^2_1\rightarrow 0,
\end{eqnarray*}
uniformly in $j$ as $n\rightarrow \infty$.
Hence, $\textcircled{3}$ holds for any weak sense limit.
\end{itemize}
\qed

\begin{remark}\label{PIn to PI}
Theorem \ref{stationarity} yields that $\{Q^{n,k}_1(\cdot)\}$ converges weakly to the probability measure of $\{U_{t_\cdot},\Pi_\cdot\}$ along the chosen subsequence.
Keeping in mind that $\rho(\cdot)$ is the unique invariant measure of $\{U_{t_\cdot},\Pi_\cdot\}$, so one can assert $\{Q^{n,k}_1(\cdot)\}$ converges weakly to $\rho(\cdot)$ along any $(n,k)$.
Formula (\ref{PIn to rho}) is therefore derived. Particularly, 
for all $\phi\in C_b({\mathcal{X}})$, we have
\[
\lim_{n\to\infty,k\to \infty}
\frac 1 k \sum^k_{j=1}\EE\langle\Pi^n_j, \phi\rangle^2
=\lim_{k\to \infty}
\frac 1 k \sum^k_{j=1}\EE\langle\Pi_j, \phi\rangle^2.
\]
\end{remark}

Finally we are ready to conclude our main result.
\begin{theorem}
For all bounded and continuous real-valued function $\phi$,
\begin{equation}
\frac 1k\sum^{k}_{j=1}\big|\langle\Pi^n_j,\phi\rangle
-\langle \Pi_j,\phi\rangle\big|^2 \to 0
\end{equation}
in probability as $n,k\to \infty$ and $n\le k$.
\end{theorem}
{\bf Proof.}
Thanks to Remark $\ref{PIn to PI}$, we have
\begin{align*}
&\lim \limits_{n\rightarrow\infty, k\rightarrow \infty}
\EE\frac 1 k \sum^k_{j=1}\big|\langle\Pi^n_j, \phi\rangle-\langle \Pi_j,\phi\rangle\big|^2\\
=& \lim\limits_{n\rightarrow\infty,k\rightarrow \infty}\Big(
\frac 1 k \sum^k_{j=1}\EE\langle\Pi^n_j, \phi\rangle^2
+\frac 1 k \sum^k_{j=1}\EE\langle\Pi_j, \phi\rangle^2
-\frac 2 k \sum^k_{j=1}\EE\langle\Pi^n_j, \phi\rangle\langle\Pi_j, \phi\rangle\Big)\\
=& \lim\limits_{k\rightarrow \infty}
\frac 2 k \sum^k_{j=1}\EE\langle\Pi_j, \phi\rangle^2
- \lim\limits_{n\rightarrow\infty,k\rightarrow \infty}
\frac 2 k \sum^k_{j=1}\EE\Big(\langle\Pi^n_j, \phi\rangle \EE\big(\phi(U_{t_j})|\mathcal{F}^Y_j\big)\Big)\\
\overset{*}{=}& \lim\limits_{k\rightarrow \infty}
\frac 2 k \sum^k_{j=1}\EE\langle\Pi_j, \phi\rangle^2
- \lim\limits_{k\rightarrow \infty}
\frac 2 k \sum^k_{j=1}\EE\langle\Pi_j, \phi\rangle\phi(U_{t_j})=0,
\end{align*}
where the equality (*) holds by formula (\ref{PIn to rho}). This implies the convergence in probability.
\qed

\section{Simulation results}
To specifically detect the source of pollution and its concentration along the river, we set $\alpha=5$, $a=1$, $b=2$, $h(z)=3z$, the observation location $x_0=0.2$, and the undesired chemical are released from location $\beta_0=0.6$ by a Poisson process with intensity $\lambda_0=10$. Let $\phi_0$ be a Dirac delta function with $\phi_0(z)=\delta(z-x_0)$, indicating that the observation data are exactly the chemical concentration at $x_0$. Suppose the initial condition $U(0,z)=1$. By Theorem $\ref{consistency of theta}$ and Newton-Raphson algorithm,  given the initial guess $\hat{\beta}_0=1/3$, $\hat{\lambda}_0=5$ and $\Delta=0.01$, we have the simulation results of $\hat{\beta}_n$, $\hat{\lambda}_n$ as Figure 1.
\begin{figure}[h]
\begin{adjustwidth}{-1cm}{}
\includegraphics[height=9cm, width=17cm]{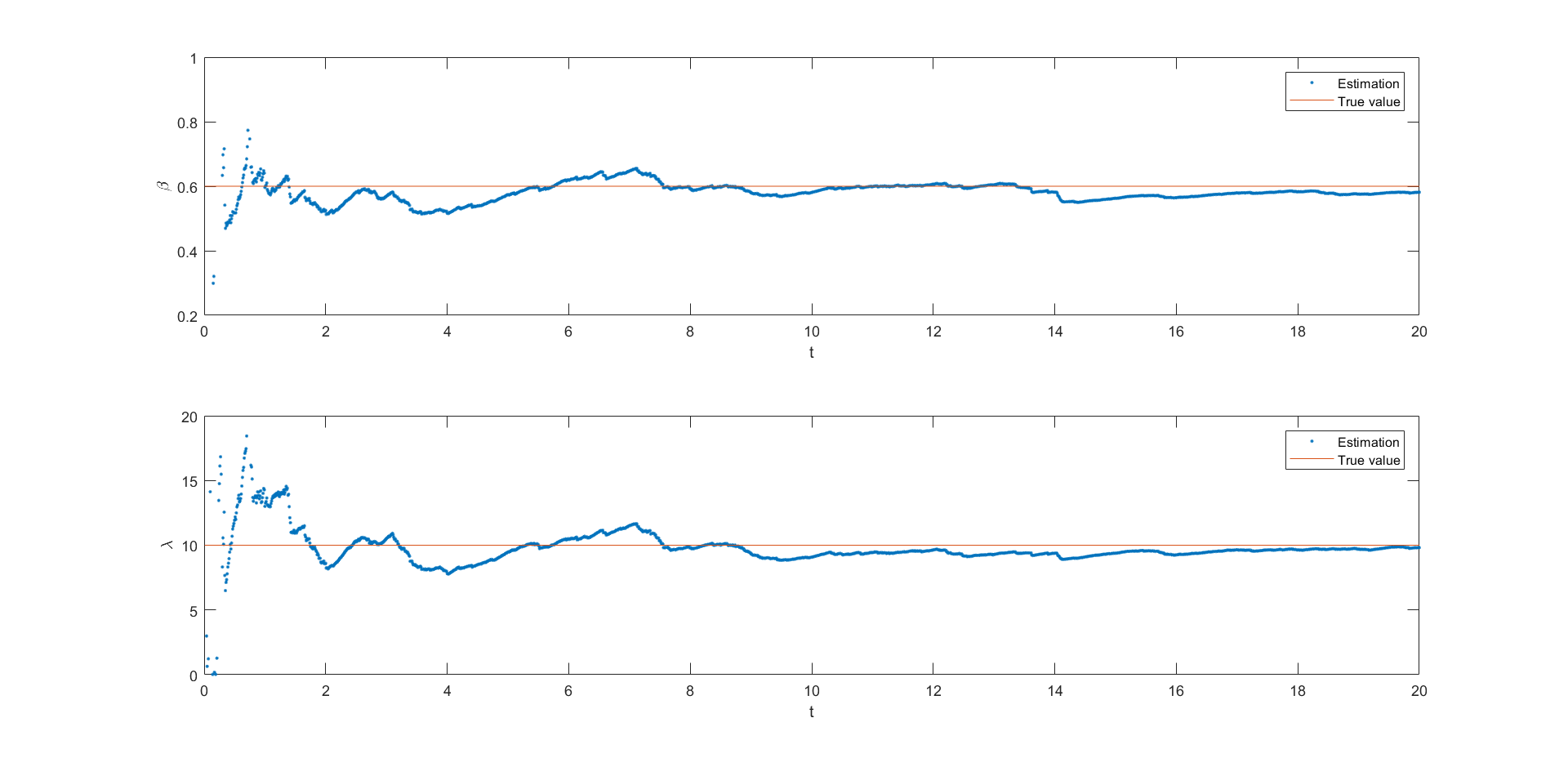}
\end{adjustwidth}\vspace{-1.2cm}
\caption{Estimation of parameters} \label{labe} 
\end{figure}

Now we are to construct the approximated filter with formula $(\ref{app fil})$. Taking $n=550$, we have the estimators based on the observation up to $t=5.5$ as $\hat{\beta}_{n}=0.597$, $\hat{\lambda}_{n}=10.038$, . Let $\phi(z)=z$, $M=500$ and define the approximated filter $\Pi^n_j$ by the sample average:
\[
\langle\Pi^n_j,\phi\rangle
=\frac {\sum^M_{l=1}\phi(\tilde U^{n,l}_{t_j})
R(\tilde U^{n,l}_{0,j},Y_{0,j})/M}
{\sum^M_{l=1}R(\tilde U^{n,l}_{0,j},Y_{0,j})/M}.
\]
where
$$
\tilde U^{n,l}_{t_j}=e^{-\alpha t_j}\int^1_0P_{t_j}(y,dx)\tilde U_0(dy)
+\int^{t_j}_0 e^{-\alpha (t_j-s)}
P_{t_j-s}(\hat\beta_n, x)d\tilde N^l_{\hat\lambda_n s},\quad l=1,\cdots,M,
$$
$$
R(\tilde U^{n,l}_{0,j},Y_{0,j})=\exp\left\{\sum^j_{i=1}[Y_{i}h(\langle \tilde U^{n,l}_{t_i},\phi_0\rangle)
-\frac 12h^2(\langle \tilde U^{n,l}_{t_i},\phi_0\rangle)]\right\}.$$
The observation data $Y_\cdot$, true concentration $U_\cdot$ and our approximation $\langle \Pi^n,\phi\rangle$ at location $x_0$ are illustrated in Figure 2.
\begin{figure}[h]
\begin{adjustwidth}{-2.1cm}{}
\includegraphics[height=11cm, width=19cm]{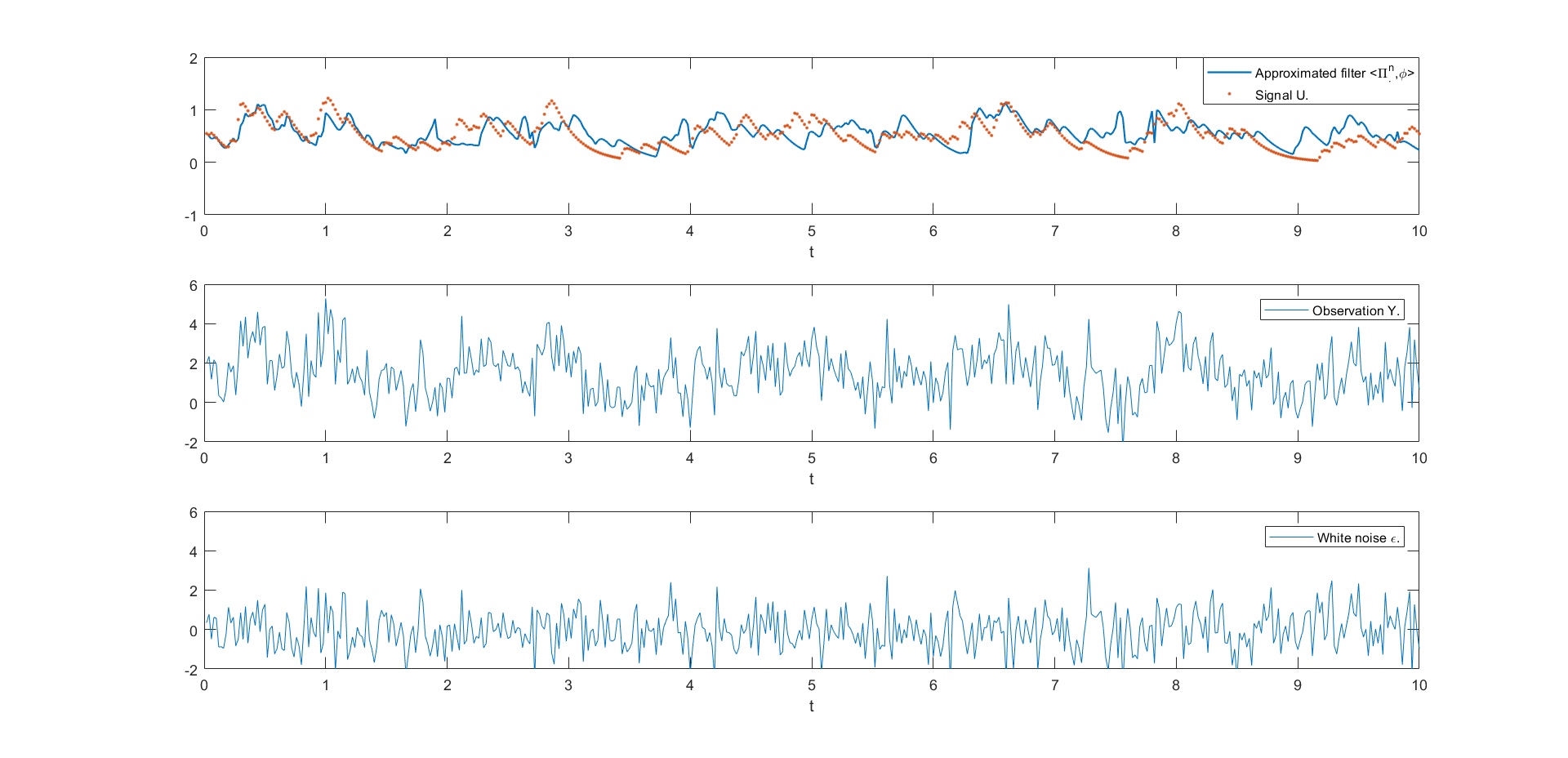}
\end{adjustwidth}\vspace{-1.2cm}
\caption{Approximation of signal process} \label{fil} 
\end{figure}

\section{Concluding remarks}
Under the background of environmental pollution detection, we put forward a nonlinear filtering problem over a long time interval where the signal process was driven by an SPDE involved with Poisson Process and unknown parameters. According to the discrete samples, the unknown parameters were figured out by the ergodic theory and statistical schemes at first. Then we proceeded to verify the uniqueness of invariant measure for the signal-filter process. Since the optimal filter is considerably hard to construct, the computed approximating filter was established. The asymptotic properties were finally argued both by theory and matlab simulations. This paper provides a comprehensive method to solve the nonlinear filtering problem. The model is of great practical value and could be widely used in many fields.

\addcontentsline{toc}{section}{Bibliography}


\begin{thebibliography}{99}

\bibitem{Bishop 2006}
Bishop, C. M. (2006) {\em Pattern Recognition and Machine Learning}. Springer-Verlag, Berlin, Heidelberg.

\bibitem{Brigo and Hanzon 1998}
Brigo, D. and Hanzon, B. (1998) On some filtering problems arising in mathematical finance. {\it Insur. Math. Econ.}, {\bf 22}, 53--64.

\bibitem{Brigo 1999}
Brigo, D., Hanzon, B. and Gland, F. L. (1999) Approximate nonlinear filtering by projection on exponential manifolds of densities. {\it Bernoulli}, {\bf 5}, 495--534.

\bibitem{Budhiraja and Kushner 1999}
Budhiraja, A. and Kushner, H. J. (1999) Approximation and limit results for nonlinear filters over an infinite time interval. {\it SIAM J. Cont. Opt.} {\bf 37}, 1946--1979.


\bibitem{Budhiraja and Kushner 2001}
Budhiraja, A. and Kushner, H. J. (2001) Monte Carlo algorithms and asymptotic problems
in nonlinear filtering. {\em Stochastics in Finite and Infinite Dimensions}, 59--87. Birkh\"auser, Boston, MA. 

\bibitem{Gordon 1993}
Gordon, N. J., Salmon, D. J. and Smith, A. F. M. (1993) Novel approach to nonlinear/non-Gaussian Bayesian state estimation. {\em IEE Proc. F}, {\bf140}, 107--113.


\bibitem{Handel 2009}
van Handel, R. (2009) The stability of conditional Markov processes and Markov chains
in random environments. {\em Ann. Probab.} {\bf 37}, 1876-1925.

\bibitem{Handel 2012}
van Handel, R. (2012) On the exchange of intersection and supremum of $\sigma$-fields in filtering theory. {\em Isr. J. Math.} {\bf 192}, 763--784.



 


\bibitem{Kallianpur and Xiong 1995}
Kallianpur, G. and Xiong, J. (1995) Stochastic Differential Equations in Infinite Dimensional Spaces.  {\em IMS Lecture Notes-Monograph Series} {\bf 26}. Institute of Mathematical Statistics, Hayward, California.


\bibitem{Kushner 1977}
Kushner, H. J. (1977) {\em Probability Methods for Approximations in Stochastic Control
and for Elliptic Equations}. Academic Press, New York.

\bibitem{Kushner and Huang 1986}
Kushner, H. J. and Huang, H. (1986) Approximation and limit results for nonlinear filters with wide bandwidth observation noise, {\em Stochastics}, {\bf 16}, 65--96.

\bibitem{Kushner 1990}
Kushner, H. J. (1990) {\em Weak Convergence Methods and Singularly Perturbed Stochastic Control and Filtering Problems}. Bonston, MA: Birkh\"auser.

\bibitem{Kushner 2008}
Kushner, H. J. (2011) Numerical approximations to optimal nonlinear filters, {\em The Oxford Handbook of Nonlinear Filtering}, D. Crisan and B. Rozovski, Eds. Oxford, U.K.: Oxford, ch. 28.

\bibitem{Mandal and Mandrekar 2000}
Mandal, P. K. and Mandrekar, V. (2000) A Bayes formula for Gaussian noise processes and its applications. {\em SIAM J. Cont. Opt.} {\bf 39}, 852--871.

\bibitem{Kouritzin 2004}
Kouritzin, M., Sun, W. and Xiong, J. (2004) Nonlinear filtering for reflecting diffusions in random enviroments via nonparametric estimation. {\em Electron. J. Probab.}, {\bf 9}, 560--574.

\bibitem{Van Leeuwen 2010}
Van Leeuwen, P. J. (2010) Nonlinear data assimilation in geosciences: An extremely efficient particle filter. {\it Q. J. Roy. Meteor. Soc.}, {\bf136}, 1991--1999.

\bibitem{Xiong 2008}
Xiong, J. (2008) {\em An introduction to stochastic filtering theory} (Vol. 18). Oxford University Press on Demand.

\end{thebibliography}
\end{document}